\documentclass{article}
\usepackage{amsfonts,hyperref,amsmath,amssymb,authblk}
\usepackage{fullpage}
\usepackage{physics}
\usepackage{enumitem}

\newcommand{\conv}[1]{\operatorname{\rm conv}\left( #1 \right)}			

\newtheorem{theorem}{Theorem}[section]

\title{A proof of the elliptical range theorem via Kippenhahn's theorem}

\author[1]{Pietro Paparella} 
\author[2]{Luis J.~Ramirez}  
\author[2]{Yen-Fen Wang}

\affil[1]{Division of Engineering and Mathematics, University of Washington Bothell, Bothell, WA 98011-8246, USA (pietrop@uw.edu)}
\affil[2]{University of Washington Bothell, Bothell, WA 98011-8246, USA}

\begin{document}
\maketitle

\begin{abstract}
The elliptical range theorem asserts that the field of values (or numerical range) of a two-by-two matrix with complex entries is an elliptical disk, the foci of which are the eigenvalues of the given matrix. Many proofs of this result are available in the literature, but most, with one exception, are computational and quite involved. In this note, it is shown that the elliptical range theorem follows from the properties of plane algebraic curves and a straightforward application of a well-known result due to Kippenhahn.
\end{abstract}



\section{Introduction}

If $A$ is an an $n$-by-$n$ matrix, then the \emph{field (of values)} or the \emph{numerical range of $A$}, denoted by \(F(A)\), is defined by \( F(A) = \left\{ x^*A x : x^*x = 1 \right\} \subseteq \mathbb{C} \). The field possesses many desirable properties and, as such, is useful in operator theory and other subjects including quantum computing (e.g., see Horn and Johnson \cite{hj1994} and references therein). In particular, the celebrated Toeplitz-Hausdorff theorem asserts that the field is convex. 

Most proofs of this result begin by observing that the problem reduces to the two-dimensional case \cite{d1971, d1957, hj1994, l1994}. In this particular case we obtain the following fundamental result, which is known as the \emph{elliptical range theorem}. 

\begin{theorem}
\label{ert}
[Elliptical range theorem] 
If A is a two-by-two matrix with complex entries and eigenvalues $\lambda_1$ and $\lambda_2$, then $F(A)$ is an elliptical disk centered at $(1/2)\trace(A)$, has foci $\lambda_1$ and $\lambda_2$, and has minor axis length equal to 
\[ \sqrt{\trace{(A^*A)} - |\lambda_1|^2 - |\lambda_2|^2}. \]
\end{theorem}

There are many proofs of this result \cite{d1957, hj1994, j1974, l1996, m1932} and most, with one exception \cite{l1996}, are computational and quite involved. In the sequel, it will be shown that the elliptical range theorem follows from properties of plane algebraic curves and a straightforward application of a well-known result due to Kippenhahn.

\section{Background}

In this section we recall relevant background.

\subsection{Algebraic Curves} 

If $p \in \mathbb{C}[x,y]$ is a polynomial of degree $n$, then the \emph{plane algebraic curve with respect to p}, denoted by $\gamma = \gamma_p$, is defined by
\[ \gamma = \left\{ (x,y) \in \mathbb{C}^2 \mid p(x,y) = 0 \right\}. \] 
If $P \in \mathbb{C}[x,y,z]$ is a homogeneous polynomial of degree $n$, then the \emph{plane projective curve with respect to P}, denoted by $\kappa = \kappa_P$, is defined by 
\[ \kappa = \left\{ (x,y,z) \in \mathbb{CP}^2 \mid P(x,y,z) = 0 \right\}. \] 
The \emph{degree of $\gamma$} (\emph{degree of $\kappa$}), denoted by $\deg \gamma$ (respectively, $\deg \kappa$), is defined by $\deg \gamma = \deg p$ (respectively, $\deg \kappa = \deg P$). 
The \emph{real part} of a plane algebraic curve $\gamma$ (plane projective curve $\kappa$) is defined by $\Re(\gamma) = \left\{ (x,y) \in \mathbb{R}^2 \mid p(x,y) = 0 \right\}$ (respectively, $\Re(\kappa) = \left\{ (x,y,z) \in \mathbb{RP}^2 \mid P(x,y,z) = 0 \right\}$).

If $p \in \mathbb{C}[x,y]$, then $H[p](x,y,z) := z^{\deg p} p(x/z,y/z)$ is a homogeneous polynomial. If $P \in \mathbb{C}[x,y,z]$ is a homogeneous polynomial, then $B[P] (x,y) := P(x,y,1)$ is a bivariate polynomial. Thus, every plane algebraic curve can be identified with a plane projective curve and vice-versa. 

If $\kappa_P$ is a plane projective curve of degree $n$, then the \emph{dual of $\kappa_P$}, denoted by $(\kappa_P)^\delta = \kappa_{P^\delta}$, is the unique plane projective curve of degree $m$ such that 
\begin{equation} P^\delta(u, v, w) = 0 \label{tangentialequation} \end{equation}
if and only if the line 
\begin{equation} ux + vy + wz = 0 \label{tangentline} \end{equation} 
is tangent to $\kappa_P$. 
It is well-known that $((\kappa_P)^\delta)^\delta = \kappa_P$. 

To find the \emph{point-equation} $P(x,y,z) = 0$ given the \emph{tangential equation} $P^\delta (u,v,w)=0$, one can eliminate the variables $u,v,w,\lambda$ from \eqref{tangentialequation}, \eqref{tangentline}, and 
\begin{equation}
\label{partials}
\frac{\partial P^\delta}{\partial u} + \lambda x = 0,~\frac{\partial P^\delta}{\partial v} + \lambda y = 0,~\frac{\partial P^\delta}{\partial w} + \lambda z = 0
\end{equation} 
(see, e.g., Salmon \cite[p.~76]{s1960}). 

\subsection{The Field of Values}

Let $A$ be an $n$-by-$n$ matrix with complex entries. The following properties are well-known and otherwise easy to establish:
\begin{enumerate}
[leftmargin=.5in, label = {\bf P\arabic*}]
\item \label{P1} $F(\alpha A + \beta I) = \alpha F(A) + \beta$ \cite[Properties 1.2.3 \& 1.2.4]{hj1994}.
\item \label{P2} If \( A \) is \emph{normal}, i.e., if \( A^*A = AA^*\), then \( F(A) = \conv{\sigma(A)} \), in which \(\sigma(A)\) denotes the spectrum of $A$ \cite[Property 1.2.9]{hj1994}.
\item \label{P3} If $U$ is a \emph{unitary matrix}, i.e., if $U^*U = UU^* = I_n$, then $F(U^* A U) = F(A)$. 
\end{enumerate}

The following result is due to Kippenhahn \cite{k1951}. 

\begin{theorem}
[{\cite[Theorem 10]{k1951}}]
If $A$ is an $n$-by-$n$ matrix with complex entries, then there is a plane algebraic curve $\kappa_P$ of class $n$ such that 
\[ F(A) = \conv{\Re(\gamma_{P(x,y,1)})}. \] 
Furthermore, if $H_1 := (A + A^*)/2$ and $H_2 := (A-A^*)/(2i)$, then 
\begin{equation}
P^\delta = \left \vert H_1 u + H_2 v + I_n w \right \vert. 
\end{equation} 
\end{theorem}     

\section{The Proof}

\emph{Proof of Theorem \ref{ert}}. Let $A$ be a two-by-two matrix with complex entries and eigenvalues $\lambda_1$ and $\lambda_2$. We follow the cases established by \cite{l1996}. 

If $A$ is normal, then, as a consequence of \ref{P2}, $F(A) = \conv{\lambda_1,\lambda_2}$, which can be viewed as an  ellipse with foci $\{ \lambda_1, \lambda_2 \}$ and minor axis length of zero.      
    
Suppose $A$ is not normal. Without loss of generality, we may assume that $\trace{A} = 0$ (othwerwise, we can replace $A$ with $A - (\trace{A}/2) I_2$ in view of \ref{P1}). By the Schur decomposition theorem, the matrix $A$ is unitarily similar to the matrix  
\begin{equation*}
B=
\begin{bmatrix}
\lambda & b \\
0 & -\lambda
\end{bmatrix},~b \ne 0.
\end{equation*}

Consider the following cases:
\begin{enumerate}
[leftmargin=.5in, label=(\roman*)]
\item $\lambda = 0$. Without loss of generality, we only consider the case when $b = 1$ (if $b \ne 1$, one may consider the matrix $B/b$ in view of \ref{P1}). A straightforward computation shows that
\begin{equation*}
P^\delta = |H_{1} u + H_{2} v + I_{n} w| = -\frac{u^2}{4} - \frac{v^2}{4} + w^{2}. 
\end{equation*}
Eliminating the variables $u$, $v$, $w$, and $\lambda$ from \eqref{tangentialequation}, \eqref{tangentline}, and \eqref{partials} yields the plane projective curve
\begin{equation*}
\kappa_P = \left\{ (x,y,z) \in \mathbb{CP}^2 \mid \frac{x^2}{\frac{1}{4}} + \frac{y^2}{\frac{1}{4}} - z^2 = 0 \right\}.
\end{equation*}
The real plane algebraic curve 
\[
\Re(\gamma_{P(x,y,1)}) = \left\{ (x,y) \in \mathbb{R}^2 \mid \frac{x^2}{\frac{1}{4}} + \frac{y^2}{\frac{1}{4}} - 1 = 0 \right\}
\]
is a circle centered at the origin with diameter 
\[ 2 \left( \frac{1}{2} \right) = 1 = \sqrt{\trace{B^*B}}. \]

\item $\lambda \ne 0$. Without loss of generality, we only consider the case when $\lambda = 1$ (if $\lambda \ne 1$, one may consider the matrix $B/\lambda$ in view of \ref{P1}). A straightforward computation shows that
\begin{equation*}
P^\delta = |H_{1} u + H_{2} v + I_{n} w| = -\left(1 + \frac{b \Bar{b}}{4} \right) u^2 - \frac{b \bar{b}}{4} v^2 + w^{2}. 
\end{equation*}
Eliminating the variables $u$, $v$, $w$, and $\lambda$ from \eqref{tangentialequation}, \eqref{tangentline}, and \eqref{partials} yields the plane projective curve
\begin{equation*}
\kappa_P = \left\{ (x,y,z) \in \mathbb{CP}^2 \mid \frac{x^2}{1 + \frac{b\bar{b}}{4}} + \frac{y^2}{\frac{b\bar{b}}{4}} - z^2 = 0 \right\}.
\end{equation*}
The plane algebraic curve 
\[
\Re(\gamma_{P(x,y,1)}) = \left\{ (x,y) \in \mathbb{R}^2 \mid \frac{x^2}{1 + \frac{b\bar{b}}{4}} + \frac{y^2}{\frac{b\bar{b}}{4}} - 1 = 0 \right\}
\]
is an ellipse centered at the origin with minor axis length equal to 
\[ 2 \left( \frac{\sqrt{b\bar{b}}}{2} \right) = b \bar{b} = \sqrt{\trace{B^*B} - 1^2 - (-1)^2} \]
and foci
\[ \pm c = \sqrt{\left(1 + \frac{b\bar{b}}{4}\right) - \frac{b\bar{b}}{4}} = \pm 1.\]
\end{enumerate}
    
\bibliographystyle{abbrv}
\bibliography{ert_kipp}

\end{document}